\RequirePackage{fix-cm}
\documentclass[smallextended]{svjour3}       
\smartqed  
\usepackage{graphicx}
\usepackage{color}
\usepackage{amsmath,amssymb}

\newtheorem{Example}{Example}
\newtheorem{Theorem}{Theorem}
\newtheorem{Lemma}{Lemma}
\newtheorem{Remark}{Remark}

\newcommand\mybox{\hbox to 0pt{}\hfill$\rlap{$\sqcap$}\sqcup$}
\begin{document}

\title{On Upper Approximations of Pareto Fronts}

\author{I. Kaliszewski \and J. Miroforidis}


\institute{I. Kaliszewski \at Systems Research Institute, Polish
Academy of Sciences, ul. Newelska 6, 01-447 Warsaw, Poland. \\
Warsaw School of Information Technology, ul. Newelska 6, 01-447 Warsaw, Poland. \\
              Tel.: +48 22 3810392 \\
              \email{ignacy.kaliszewski@ibspan.waw.pl}  \\
              Corresponding author
           \and
           J. Miroforidis \at Systems Research Institute, Polish
Academy of Sciences, ul. Newelska 6, 01-447 Warsaw, Poland. \\
}

\date{Received: date / Accepted: date}

\maketitle

\begin{abstract}

In one of our earlier works, we proposed to approximate Pareto fronts to multiobjective optimization problems by two-sided approximations, one from inside and another from outside of the feasible objective set, called, respectively, lower shell and upper shell. We worked there under the assumption that for a given problem an upper shell exists. As it is not always the case, in this paper we give some sufficient conditions for the existence of upper shells.

We also investigate how to constructively search infeasible sets to derive upper shells. We approach this issue by means of problem relaxations. We formally show that under certain conditions some subsets of lower shells to relaxed multiobjective optimization problems are upper shells in the respective unrelaxed problems.

Results are illustrated by a numerical example representing a~small but real mechanical problem. Practical implications of the results are discussed.

\end{abstract}

\keywords{Multiobjective optimization \and Global Pareto optimum \and Two-sided Pareto front approximations \and Pareto front upper approximation existence \and Evolutionary multiobjective optimization}

\section{Introduction}
\label{section1}

In Multiobjective Optimization (MO) problems, any finite subset of the feasible set is a lower (we mean: feasible) discrete approximation of the efficient set. Such subsets, if containing no dominated (in the sense of Pareto) elements, are called \textit{lower shells}.

When coupled with dual constructs, namely \textit{upper shells}, lower and upper shells produce (via objective function mappings) the so-called two-sided Pareto front approximations \cite{KM_2014}.
In inexact MO, two-sided Pareto front approximations enable a natural, qualitative definition of global Pareto optimum: \textit{the global Pareto optimum is reached if a two-sided Pareto front approximation is known with its lower and upper part close enough (various metrics has been proposed) to each other}~\footnote{\mbox{ }By an analogy to singleobjective optimization, where ''local optimum'' is not necessarily the true (global) optimum, the majority of inexact MO methods solve MO problems ''locally'', i.e., not necessarily solutions they provide are true Pareto optima (true Pareto fronts). This as a consequence of method inexactness, but such solutions can be regarded as solutions to the problem only if they satisfy the above definition of global Pareto optimum.}. This definition can be quantified as MO problem contexts dictate.

Roughly speaking, an \textit{upper shell} (a finite set) approximates the efficient set from outside of
the feasible set. It is required that no element of the upper shell is dominated by any element of the efficient set, a natural prerequisite for any Pareto front (i.e., the image of the efficient set under the objective function mapping) approximation of that kind. But this means that verification of that requirement necessitates the knowledge of the efficient set. Hence, the definition of upper shell is not, in general, operational.

A weaker notion is \textit{upper approximation} (a finite set). An upper approximation is defined with respect to a given lower shell. The aforementioned requirement for upper shells (no element of the upper shell is dominated by any element of the efficient set) is weakened to the requirement that no element of the upper approximation is dominated by any element of that given lower shell. In consequence, the definition of upper approximation becomes operational.

However, such a weakening may cause that some elements of an upper approximation are dominated by some elements of the efficient set, definitely a~harmful property. It is of interest then to identify instances in which upper approximations coincide with upper shells and such instances have been identified in this paper.

In our earlier works \cite{KM_2014,KMP_2012},
where we were concerned with algorithmic issues of the derivation of lower shells and upper shells, we worked under the assumption that for a given problem an upper shell exists. However, an upper shell not always exists (see example in Section \ref{section7}). A~number of problems where no upper shell exists have been identified \cite{KM_2010,KM_2012}.
Therefore, the problem of existence of upper shells deserves consideration. In this work, we give sufficient conditions for an upper shell to exist.

Another issue is how to search the infeasible set to derive upper shells. We approach that issue by means of \textit{problem relaxations}. We show that some subsets of lower shells in relaxed MO problems are upper shells in the respective unrelaxed problems.

The practical importance of lower and upper shells lies in that they enable calculation of lower and upper bounds on values of objective functions for any implicit efficient solution (an implicit efficient solution is an efficient solution which can be derived by solving a scalarized MO problem, but as long as the problem is not solved, it remains unknown) \cite{KMP_2012,KMP_2011}.
In calculations of upper bounds on objective function values, upper shells can be replaced by upper approximations if it is known that they coincide. This paves the way for approximate (inexact) computations with controlled accuracy in MO, with applications to Multiple Criteria Decision Making \cite{KMP_2012,KMP_2011}.
This fact is the main motivation of our research presented here -- we would like to identify conditions under which upper bounds on values of objective functions can be calculated (any lower shell enables calculation of lower bounds). Thus, our results presented below are of existential type. Considerations relating to quality issues of two-sided Pareto front approximations have been discussed in \cite{KM_2014,KMP_2012}.

The outline of the paper is as follows. Section \ref{section2} relates our work to the relevant literature of the field.  In Section \ref{section3}, we present preliminaries. In Section \ref{section4}, we identify instances in which upper approximations coincide with upper shells. In Section \ref{section5}, we show how upper approximations, whenever they  exist, can be derived by relaxations of  MO problems, whereas in Section \ref{section6} we discuss invariance of upper approximations under order invariant transformations of objectives. In Section
\ref{section7}, we give sufficient conditions for existence of upper shells, and in Section \ref{section8}, we present an application of these results to a practical problem. Section \ref{section9} concludes.

\section{Related works}
\label{section2}

The need for Pareto front (PF) approximations was realized in the MO community early. Following
\cite{RW_2005}, approximation concepts can be divided into those based on exact MO methods and inexact ones (i.e., all kinds of heuristics, evolutionary computations including). According to the classification scheme proposed in \cite{RW_2005}, $0$th order approximations are discrete (pointwise) approximations and consist of a number (usually limited) of elements of PFs generated by a solution method. Higher order approximations consist of some constructs built on those elements, cf. e.g. \cite{KMW_2003}.

The survey \cite{RW_2005} concentrates on exact methods and covers the period of 1975-2005. Since 2005, other approximation concepts based on exact methods have been proposed \cite{LBKM_2005,BLL_2014,Lotov_2015,HMW_2012,Evtushenko_Posypkin_2013,Bradley_2014}. In \cite{LBKM_2005,BLL_2014}, the efficient set is approximated from inside of the feasible set by solving series of optimization problems. The approach has been recently refined in \cite{Lotov_2015} by an application of the decomposition principle. The method proposed in \cite{HMW_2011,HMW_2012} interpolates a number of elements of the Pareto front, and this interpolation gives rise to a mixed integer linear surrogate problem.
Properties of $\varepsilon$-Pareto set approximations, once they are given, are investigated in \cite{Evtushenko_Posypkin_2013}. In the same work, a method to derive $\varepsilon$-Pareto set approximations with the Lipschitz type information, extending earlier works in that direction (\cite{Evtushenko_Potapov_1987}) is given. A statistical model (Kriging) is applied in \cite{Bradley_2014} to facilitate feasible solution set sampling in a quest for the Pareto front. All those works provide higher than $0$th order approximations. Of higher that $0$th order approximations are also those which hybridize exact and inexact approaches (\cite{LBKM_2005,BLL_2014}). All those methods also explicitly or implicitly assume (with an exception for \cite{Evtushenko_Posypkin_2013}) that the feasible set has interior, thus excluding combinatorial problems from considerations.

However, works on PF approximations in discrete (combinatorial) problems are also represented in the literature. In \cite{V_K_1981}, definitions of lower and upper bound sets, which coincide with the definition of lower shell and is quite close to the definition of upper shell used in this work, is introduced in the context of multiobjective knapsack problem with integer variables, solved via dynamic programming. Approximations of the Pareto front for the general case (no assumption on the problem considered), based on deriving efficient elements by solving MO problems scalarized by the Chebyshev function, are proposed in \cite{Kaliszewski_2006}. In \cite{EG_2007}, the idea of lower and upper bound sets is applied to multiobjective combinatorial problems; bound sets are derived by solving a~number of scalarized problems. The same approach as in \cite{Evtushenko_Posypkin_2013,Evtushenko_Potapov_1987}
is applied to mixed integer nonlinear problems in
\cite{Evtushenko_Posypkin_2011}. Recently, an idea of cover sets, based on the dominance relation, to represent PFs has been elaborated in \cite{VWW_2016}.

In \cite{Kaliszewski_2016}, lower and upper shells were applied to provide bounds on optimal solutions to biobjective knapsack problems in cases commercial mixed-integer programming solvers (like CPLEX) hit time or memory limits.

In all those works, with no exception, no attempt is made to exploit information which is provided by some specific infeasible solutions. In contrast to that, in this work we follow the other course, namely we are interested in approximations based exclusively on inexact MO methods (to ensure generality of the course, no assumption is made on whether the feasible set has interior). We have been inspired by the success story of population based methods \cite{Deb_2001,Coello_2002} (in the MO domain customarily termed Evolutionary Multiobjective Optimization -- EMO), when applied to a wide range of practical problems, cf. e.g. \cite{Coello_2004,Talbi_2009,Di_Barba_2010,GHMC_2013,QWPS_2013,RLRS_2015}. Population based methods, though by their nature inexact, have gained much popularity in application oriented communities which have no problem with accepting suboptimal solutions in exchange for method generality, versatility and simplicity, allowing easy in-house codings. 

Following the classification given in \cite{RW_2005}, one can perceive EMO methods as population based $0$th order approximations, however with the distinction that EMO methods, as a rule, produce lower (we mean: feasible) approximations of PFs with no guarantee that they include any PF elements. By this, in contrast to approximations based on exact methods, EMO methods have no built-in ''secure anchors'' in PFs. Without knowing, at least some, elements of the PF, there is no trustworthy measure of accuracy of EMO approximations. This fact inspired the authors to investigate the possibility to provide two-sided PF approximations which give rise to such measures (\cite{KM_2014,KMP_2012,KM_2010,KM_2012}). The authors are aware of only one work in which a similar reasoning is present \cite{Legriel_et_all_2010}; however, in that work two-sided approximations were not generated intentionally, as the authors pursuit in their works.


Another reason why we focus on approximations based on inexact MO methods is that such methods seem to rise less concerns to the question of scalability than their exact MO method based counterparts. In population based approximations,  we have no formal constructs to recalculate/update, except objective functions or a fitness function built on them. So it seems that in large multiobjective optimization problems  the “curse of dimensionality” will trouble the exact MO methods and exact method based approximations to a~much larger extent than the population based approximations.

\section{Preliminaries}
\label{section3}

Consider the MO problem
\begin{equation}
\label{eq3.1}
\begin{array}{c}
 ''max'' f(x) \\ x \in X_0 \, ,
 \end{array}
\end{equation}
where $X_0 \subset {\cal R}^n$ is a compact (i.e., closed and bounded) set, $f: {\cal R}^n \rightarrow {\cal R}^k, \
f=(f_1,\dots,f_k), \ f_l: {\cal R}^n \rightarrow {\cal R}, \
l=1,\dots,k, \ k \geq 2$, $f_l$ are objective functions; $"max"$
denotes the operator of deriving the set (denoted $N$) of  \textit{efficient} (as
defined below) elements of $X_0$.
We assume that $N$ is not empty.

Below we will use the following notation: $P = f(N), \ Z = f(X_0)$ (feasible objective set), $R^k_+ = \{ y \in {\cal R}^k \,|\, y_l \geq 0, \
l=1,\dots,k \}$. Set $P$ is called \textit{Pareto front}.

\textit{Dominance relation} $\prec$ on ${\cal R}^n$ is defined as
\[
 x \prec x' \ \Leftrightarrow \ f(x) \ll f(x'),
 \] where $\ll$
denotes $f_l(x) \leq f_l(x'), \ l = 1,\dots,k$, and $f_l(x) <
f_l(x')$ for at least one $l$.

Elements $x$ of $X_0$ for which there exists no $x' \in X_0$ such that $x \prec x'$, are called \textit{efficient}.


In subsequent sections, we will refer to the concept of \textit{lower shell} and \textit{upper shell}, already presented in a series of publications (\cite{KM_2014,KMP_2012,KM_2010,KM_2012}).

\vspace{0.5cm}
\textit{Lower shell} is a finite nonempty set $S_L \subseteq
X_0$, elements of which satisfy
\begin{equation}
\label{eq3.2}
\forall \, x \in S_L \ \  \not\exists \, x' \in S_L \ \ x
\prec x' \,.
\end{equation}

The condition (\ref{eq3.2}) ensures that a lower shell does not contain redundant (in the sense of the dominance relation) elements. No element of $S_L$ is dominated by any other element of this set. In terms of relations, set $S_L$ consists only of elements which are maximal to relation $\prec$\,. In [11], sets satisfying (\ref{eq3.2}) have been called \textit{inherently nondominated}.

Element $y^{nad}$ is defined as \[ y^{nad}_l =
\min_{x \in N} f_l(x), \ l=1,\dots,k \,.
\]

\textit{Upper shell} is a finite nonempty set $S_U
\subseteq {\cal R}^n \setminus X_0$, elements of which
satisfy
\begin{equation}
\label{eq3.3}
\forall \, x \in S_U \ \  \not\exists \, x' \in S_U \ \ x' \prec
x \,,
\end{equation}
\begin{equation}
\label{eq3.4}
\forall \, x \in S_U \ \  \not\exists \, x' \in N \ \ x
\prec x' \,,
\end{equation}
 \begin{equation}
 \label{eq3.5}
\forall \, x \in S_U \ \  y^{nad} \ll f(x) \,.
\end{equation}

Condition (\ref{eq3.3}) ensures that an upper shell does not contain redundant (in the sense of the dominance relation) elements. No element of $S_U$ dominates any other element of this set.

Condition (\ref{eq3.4}) ensures that no element of an upper shell is dominated by an element of set $N$ (i.e., by an efficient element).

Condition (\ref{eq3.5}) precludes inclusion into upper shells elements which in no circumstances can dominate an element of $N$.

Element $y^{nad}(S_L)$ is defined as \[ y^{nad}_l(S_L) =
\min_{x \in S_L} f_l(x), \ l=1,\dots,k \,.
\]

Given a lower shell $S_L$, \textit{upper approximation} is a finite nonempty set $A_U
\subseteq {\cal R}^n \setminus X_0$, elements of which
satisfy

\begin{equation}
\label{eq3.31}
\forall \, x \in A_U \ \  \not\exists \, x' \in A_U \ \ x' \prec
x \,,
\end{equation}
\begin{equation}
\label{eq3.41}
\forall \, x \in A_U \ \  \not\exists \, x' \in S_L \ \ x \prec x' \,, 
\end{equation}
 \begin{equation}
 \label{eq3.51}
\forall \, x \in A_U \ \  y^{nad}(S_L) \ll f(x) \,.
\end{equation}

Condition (\ref{eq3.31}) plays the same role as condition (\ref{eq3.3}).

Condition (\ref{eq3.41}) and condition (\ref{eq3.51}) are consequences of the fact that in general set $N$ is not known.

\begin{Remark}
\label{remark3.1}
An upper shell is an upper approximation (the opposite statement does not hold). Hence, a problem with no upper approximation possesses no upper shell.
\end{Remark}

Intuitively, upper approximations are meaningful only if $S_L + \varepsilon \cong N$ with $\varepsilon > 0$ sufficiently small, but such an intuition is valid only for sets with interior. In the case of discrete sets, one should rather work with upper shells.

\section{MO Instances where Upper Approximations are Upper Shells}
\label{section4}

Lower shells are meant to be one-sided representations of $N$ from inside of the feasible set.
Similarly, upper shells and upper approximations are meant to be one-sided representations of $N$ from outside of the feasible set.

The following lemmas identify cases where condition (\ref{eq3.41}) implies condition (\ref{eq3.4}), hence upper approximation $A_U$ is an upper shell $S_U$. No assumption about the nature of the underlying problem, such as continuity, discreetness, convexity or connectivity, is made here.

\vspace{0.2cm}
\begin{Remark}
\label{remark4.1}
The dominance relation $\bar{x} \prec \tilde{x}$ holds $\Leftrightarrow$ $f(\bar{x}) \not= f(\tilde{x})$ and $ f(\bar{x}) \in f(\tilde{x})  - R^k_+$\,.
\end{Remark}

\vspace{0.2cm}
\begin{Lemma}
\label{lemma4.1}
An upper approximation $A_U$ is an upper shell only if
\[
A_U \subseteq \{ x \ | \ f(x) \in int({\cal R}^k \setminus (P - R^k_+)) \}  \,.
\]
\end{Lemma}

\noindent Proof. The proof follows as an immediate consequence of the definition of upper shell and Remark \ref{remark4.1}.

\vspace{0.2cm}
\begin{Remark}
\label{remark4.2}
Sets
\[
\{ x \ | \ f(x) \in int({\cal R}^k \setminus (P - R^k_+)) \}
\]
and
\[
\{ x \ | \ f(x) \in (P - R^k_+) \}
\]
are disjoint. 
\end{Remark}


\setcounter{figure}{0}

\vspace{0.2cm}
Figure 1 and Figure 2 give a graphical interpretation of Lemma \ref{lemma4.1}.

\begin{figure}[h]
\label{figure1}
\begin{center}
\rotatebox{0}{\includegraphics[height=2.0in]{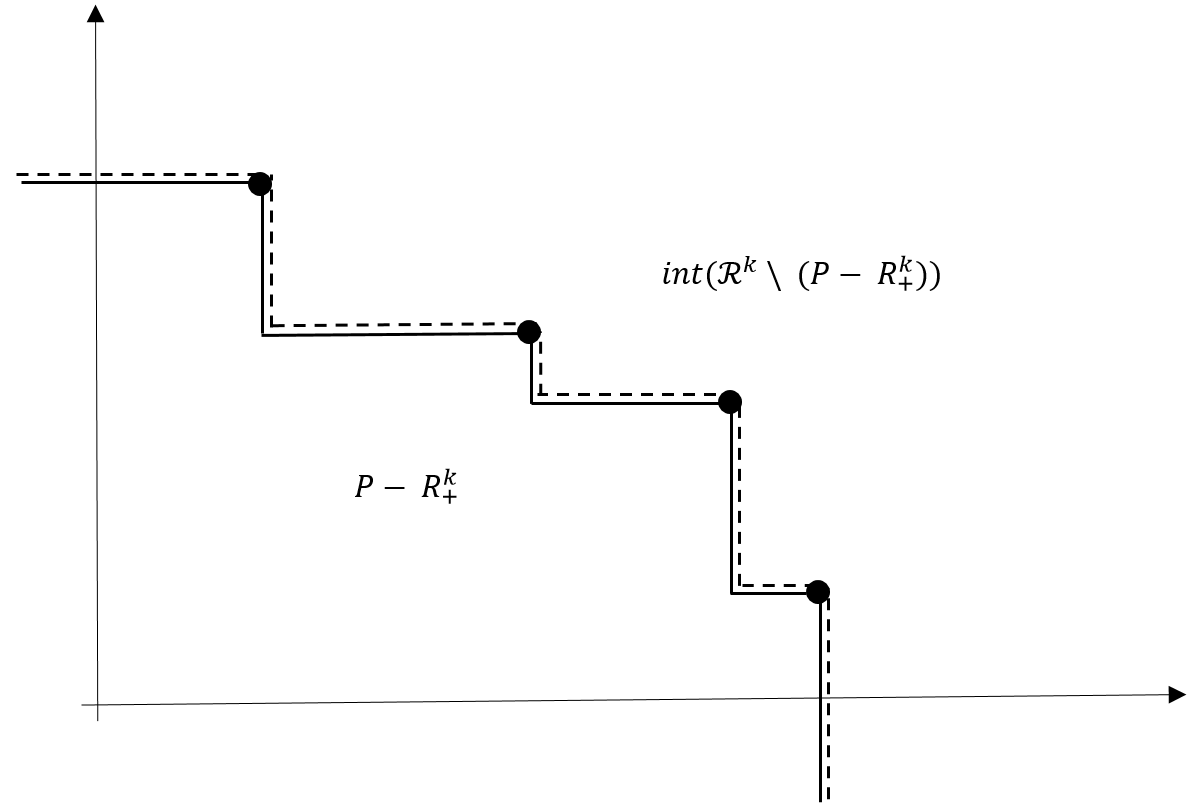}}
\end{center}
\caption{An illustration to Lemma \ref{lemma4.1}, the discrete case; bullets -- elements of $P$.}
\end{figure}

\begin{figure}[h]
\label{figure2}
\begin{center}
\rotatebox{0}{\includegraphics[height=2.0in]{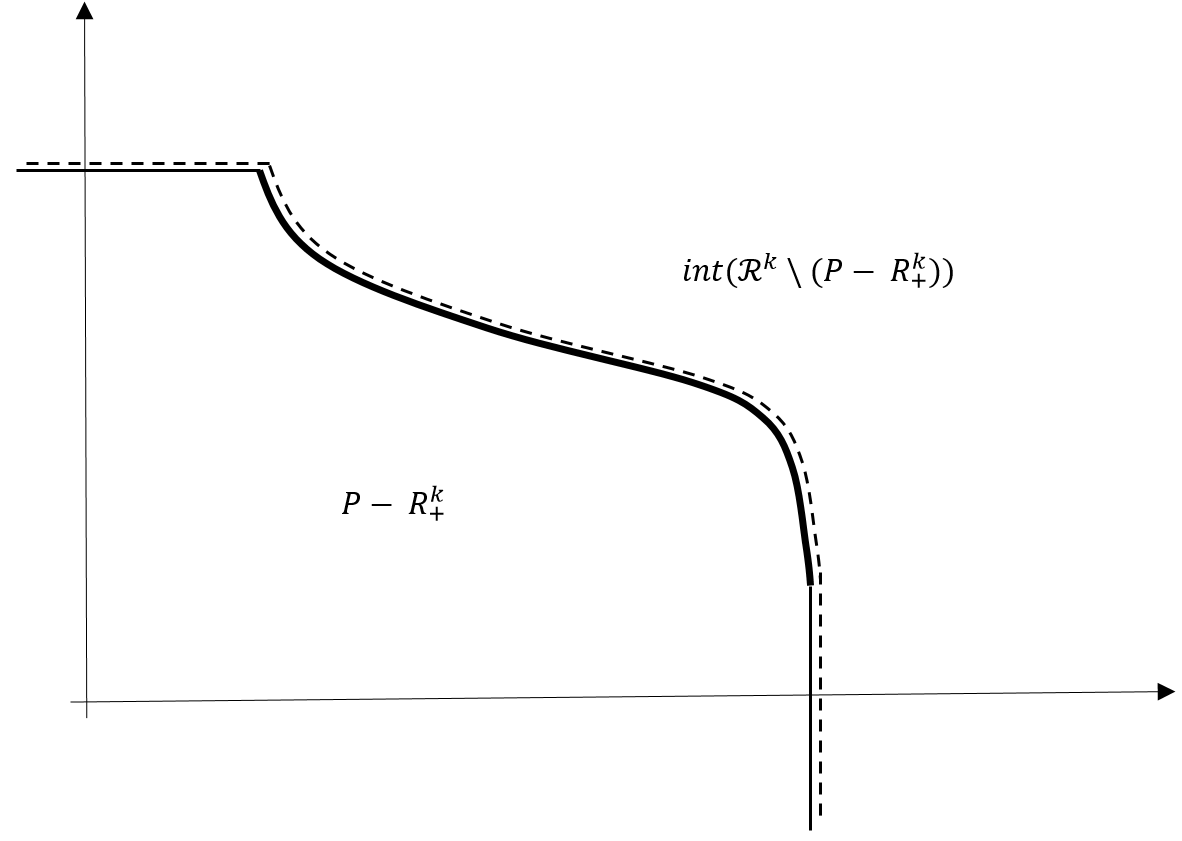}}
\end{center}
\caption{An illustration to Lemma \ref{lemma4.1}, the continuous case; thick line -- $P$.}
\end{figure}

\begin{Lemma}
\label{lemma4.2}
An upper approximation is an upper shell 
only if
\[
Z \cap  \{ f(x) \ | \ f(x) \in f(A_U) + R^k_+ \}  = \emptyset \,.
\]
\end{Lemma}

\noindent Proof. Since $Z \subseteq P - R^k_+$, the assertion of the lemma follows from Lemma \ref{lemma4.1} and Remark \ref{remark4.2}.
\mybox

\vspace{0.2cm}
\begin{Lemma}
\label{lemma4.3}
Any subset of the inverse image of any set in $P + int(R^k_+) $, \linebreak elements of which satisfy condition (\ref{eq3.3}) and condition (\ref{eq3.5}), is a valid upper shell.
\end{Lemma}

\noindent Proof. The proof follows immediately from the definition of upper shell.

\vspace{0.2cm}
Lemma 3 relates the concept of upper shells and upper approximations to works on higher than $0$th order PF approximations mentioned in Introduction, with hyperplanes as the simplest construct. Indeed, any construct  in \linebreak $P+int(R^k_+) $ which satisfies Lemma \ref{lemma4.3} is an upper shell. However, in this work we are concerned with $0$th order (pointwise) PF approximations. 

Lemma \ref{lemma4.1}--\ref{lemma4.3} give no specific guidelines how to select $A_U$ from ${\cal R}^n \setminus X_0$ to satisfy the respective assumptions for $A_U = S_U$ to hold. In general, this issue has to be investigated for each MO problem individually. However, there are classes of problems with the property that any subset of the infeasible set containing no dominating elements is an upper shell. We discuss this in Section \ref{section7}.

Now we turn to the problem of derivation of upper approximations.

\section{Derivation of Upper Approximations}
\label{section5}

Usually it is not known which part of ${\cal R}^n \setminus X_0$ should be searched for elements of $A_U$. However, upper approximations can be derived from some relaxations of problem (1), as shown by Lemma \ref{lemma5.1}.

Let $X_0 \subset X_0' \subseteq {\cal R}^n$. Consider the problem
\begin{equation}
\label{eq5.1}
\begin{array}{c} ''max'' f(x) \\ x \in X_0' \,,
\end{array}
\end{equation}
with the set of efficient elements $N'$.

Let $S_L$ denote a lower shell for problem (\ref{eq3.1}) and $S_L'$ denote a lower shell
for problem (\ref{eq5.1}).

\vspace{0.5cm}
\begin{Lemma}
\label{lemma5.1}
Let $\Theta \subseteq S_L'$ and let
elements $x$ of $\Theta$ satisfy
\begin{equation}
\label{eq5.2}
x \not\in X_0 \,,
\end{equation}
\begin{equation}
\label{eq5.3}
\forall x \in \Theta \ \ \not\exists x' \in S_L \ \ x \prec
x' \,,
\end{equation}
\begin{equation}
\label{eq5.4}
y^{nad}(S_L) \ll f(x) \,.
\end{equation}
Then $\Theta$ is an upper approximation for problem (\ref{eq3.1}).
\end{Lemma}

\noindent Proof.
By condition (\ref{eq5.2}), $\Theta \subseteq {\cal R}^n \setminus X_0$\,.

Moreover, since $\Theta$ is a subset of $S_L'$\,, by the definition of lower shell (formula (\ref{eq3.2}))
\[ \forall x
\in \Theta \ \ \not\exists x' \in \Theta \ \ x' \prec x \,,
\]
i.e., $\Theta$ satisfies condition (\ref{eq3.31}).

By condition (\ref{eq5.3}), $\Theta$ satisfies condition (\ref{eq3.41}). By
condition (\ref{eq5.4}), $\Theta$ satisfies condition (\ref{eq3.51}).
Hence, $\Theta$ satisfies the definition of upper approximation.
\mybox

\section{Invariance of approximations}
\label{section6}

Considerations of this section apply to upper approximations as well as to lower shells.

Let us consider problem (\ref{eq3.1}) with $f(x)$ replaced by some $f'(x)$, i.e., the problem

\begin{equation}
\label{eq6.1}
\begin{array}{c}
''max'' f'(x) \\ x \in X_0 \, .
 \end{array}
\end{equation}


We say that lower shells (upper approximations) are \textit{problem invariant} if every lower shell (upper approximation) to problem (\ref{eq3.1}) is a lower shell (an upper approximation) to problem (\ref{eq6.1}).

\begin{Lemma}
\label{lemma6.1}
Let for each $l=1,\dots,k$, $f_l(x)$ and $f'_l(x)$ generate the same linear order on ${\cal R}$.
Then, lower shells and upper approximations are problem invariant.
\end{Lemma}

\noindent Proof.
The proof follows from the fact that since for each $l=1,\dots,k$, $f_l(x)$ and $f'_l(x)$ generate the same linear order on ${\cal R}$, $f(x)$ and $f'(x)$ generate the same dominance relation. Since both problems have the same set $X_0$\,, every lower shell in problem (\ref{eq3.1}) satisfies the definition of lower shell in problem (\ref{eq6.1}). The same is true for the upper approximation part of the lemma. \mybox

By Lemma \ref{lemma6.1}, under specified conditions two-sided approximations in problem (\ref{eq3.1}), which could have been derived with a significant computational effort, are with no change two-sided approximations in problem (\ref{eq6.1}).

To stress the importance of the concept of invariance of approximations, it is worth mentioning that in a similar manner invariance of the efficient set ($N$) was investigated in \cite{Dumitru_Luban_1986}
 and successfully exploited in multiobjective optimization of radiotherapy planning in \cite{Romeijn_2004}.
From the latter work we learn that out of many functions proposed to measure the so-called \textit{tumor control probability} in organs to be protected against excessive radiation, one has the form as follows:
\begin{equation}
\label{eq6.2}
gEUD(d) = ( \frac{1}{v}\sum^v_{j=1} d^a_j )^{\frac{1}{a}} \,,
\end{equation}
where $d_j$ is radiation dose deposited in voxel (an element of a 3D mesh) $j$, $v$ is the number of voxels in the protected organ, $a$ is a parameter, $1 \leq a \leq \infty$. Since all $d_j$ are nonnegative, as physics dictates, clearly this function and the linear function
\begin{equation}
\label{eq6.3}
gEUD(d) = ( \frac{1}{v}\sum^v_{j=1} d_j) \,,
\end{equation}
produce, with other objective functions unchanged, the same efficient set $N$. By the same argument (Lemma \ref{lemma6.1}), both functions produce, with other objective functions unchanged, invariant lower shells and invariant upper approximations. Since in optimization problems related to oncological radiology the number of voxels depend on the mesh resolution and can reach hundreds of thousands, a simple function replacement can result in significant savings in computation load\footnote{\mbox{ }The average speed-up when calculating the value of function (\ref{eq6.3}) instead of function (\ref{eq6.2}) (averaged over $20\,000$ calculations, $a=3$) on an off-the-shelf laptop ranges $23$.}.

\section{Identification of Problems with Upper Shells}
\label{section7}

\vspace{0.5cm}
In this section, we investigate existence of upper shells. As mentioned already, there exist problems without upper shells. This fact is illustrated by the following example.

\begin{Example}
\label{ex7.1}
Let us consider the following problem
\begin{equation}
\label{eq7.1}
\begin{array}{c}
''max'' f(x) =  \left(
                            \begin{array}{c}
                            f_1(x) = -(x_1 - 3)^2 - (x_2 - 4)^2 \\ \\
                            f_2(x) = -(x_1 - 4)^2 - (x_2 - 1)^2
                            \end{array}
                       \right) \\ \\
                  X_0 =  \{\,x \ | \ 1 \leq x_1 \leq 5, \ 1 \leq x_2 \leq 5 \, \} \,,

\end{array}
\end{equation}

and its relaxation
\begin{equation}
\label{eq7.2}
\begin{array}{c}
''max'' f(x) =  \left(
                            \begin{array}{c}
                            f_1(x) = -(x_1 - 3)^2 - (x_2 - 4)^2 \\ \\
                            f_2(x) = -(x_1 - 4)^2 - (x_2 - 1)^2
                            \end{array}
                       \right) \\ \\
                   X'_0 = \{\, x \ | \ a \leq x_1 \leq b, \ a \leq x_2 \leq b, \ \ a < 1, \ b > 5 \, \} \,.
\end{array}
\end{equation}
The Pareto fronts of problem (\ref{eq7.1}) and all its relaxations are represented in Figure 3 and they are all the same. 
With the  maximum for $f_1(x)$ at $x = (3, 4)$  where $f(3,4) = (0, -10)$, and the maximum for $f_2(x)$ at $x = (4, 1)$ where $f(4,1) = (-10, 0)$, the only region of ${\cal R}^n$ where one function increases and the other decreases is defined by $\{ x \, | \, f_1(x) \geq -10, \ f_2(x) \leq 0, \ f_1(x) \leq 0, \ f_2(x) \geq -10 \} \subseteq X_0$. Thus, the problem (\ref{eq7.1}) has no upper shell.
\end{Example}

However, problem (\ref{eq7.1}) can be hardly regarded a constrained problem; its two objective functions attain their maxima inside the feasible set. Usually, problems which emerge from applications are constrained by a sort of budget constraint(s), witnessing limited resources, monetary or physical, and precluding objective functions attaining optima inside the feasible sets. In that sense, problem (\ref{eq7.1}) becomes a constrained problem for e.g. $X_0 =  \{\,x \ | \ 1.5 \leq x_1 \leq 2.5, \ 1.5 \leq x_2 \leq 2.5 \, \}$.


\begin{figure}[t]
\label{figure3}
\begin{center}
\rotatebox{0}{\includegraphics[height=1.5in]{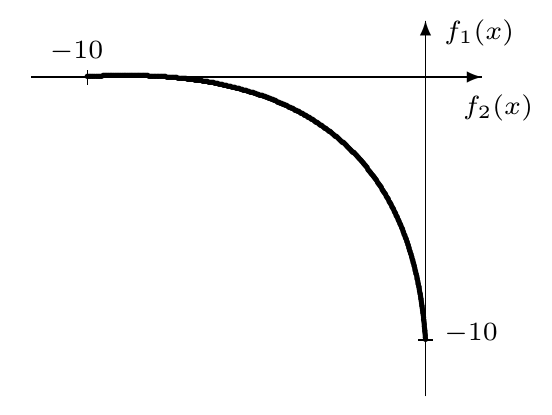}}
\end{center}
\caption{An example of problem with no upper shell. With $X_0, X'_0$ specified as in Example \ref{ex7.1}, the curve represents $f(N)$ for $X_0$ and at the same time $f(N')$ for any $X'_0$ defined as in (\ref{eq7.2}).}
\end{figure}

In general, identification of problems having upper shells is far from being trivial.
However, in some instances the existence of upper shells is relatively simple to ascertain. We recall that a function $\varphi : {\cal R}^n \rightarrow {\cal R}$, is called \textit{strongly monotonically increasing on} ${\cal R}^n$ if $x \leq x' , \ x' \not= x$, implies $\varphi(x) < \varphi(x')$ (\cite{Jahn_1986}).

Let us observe that condition (\ref{eq3.4}) is equivalent to

\begin{equation}
\label{eq7.3}
\not\exists \, x \in S_U \ \  \exists \, x' \in N \ \ x \prec x' \,.
\end{equation}

\begin{Lemma}
\label{lemma7.0}
Let objective function $f_l^*, \ l^* \in  \{1,\dots,k\}$, be strongly monotonically increasing on ${\cal R}^n$. Then, $x' \in X_0$ and $x' \leq x$, implies $x \not\prec x'$. 
\end{Lemma}

\noindent Proof. For any element $x'$ of $X_0$ and any element $x$ such that $x' \leq x$, and for strongly monotonically increasing function $f_l^*$, we have $f_l(x') < f_l(x)$.
Thus, $x \not\prec x'$. 

\begin{Lemma}
\label{lemma7.1}
Let objective function $f_l^*, \ l^* \in  \{1,\dots,k\}$, be strongly monotonically increasing on ${\cal R}^n$\,. Then, $x' \in N$, $x \in {\cal R}^n \setminus X_0$ and $x' \leq x$, implies $x \not\prec x'$. 
\end{Lemma}

\noindent Proof. Since $N \subseteq X_0$, by Lemma \ref{lemma7.0}, $x \not\prec x'$\,. Elements $x$ such that $x \in {\cal R}^n \setminus X_0$ and $x' \leq x$ exist since $X_0$ is compact.  \mybox

\vspace{0.2cm}
Lemma \ref{lemma7.1} shows how to select condidates for upper shells, which satisfy condition (\ref{eq7.3}), or equivalently, condition (\ref{eq3.4}). 
Below, we shall make use of a stronger condition, namely
\begin{equation}
\label{eq7.4}
\forall \, x \in S_U \ \  \exists \, x' \in N \ \ x' \prec x \,.
\end{equation}
Elements of $S_U$ which satisfy the above condition satisfy also condition (\ref{eq3.5}).

\begin{Lemma}
\label{lemma7.2}
Let all objective functions $f_l, \ l = 1,\dots,k$, be strongly monotonically increasing on ${\cal R}^n$\,. Then, $x' \in X_0$ and $x' \leq x$, implies $x' \prec x$.
\end{Lemma}

\vspace{0.2cm}
\noindent Proof. The proof is an immediate consequence of the assertion that all objective functions are strongly monotonically increasing on ${\cal R}^n$\,. 

\begin{Lemma}
\label{lemma7.3}
Let all objective functions $f_l, \ l = 1,\dots,k$, be strongly monotonically increasing on ${\cal R}^n$. Then, $x' \in N$ and $x' \leq x$, implies $x' \prec x$ and $x \in {\cal R}^n \setminus X_0$.
\end{Lemma}

\noindent Proof. Since $N \subseteq X_0$, by Lemma \ref{lemma7.2}, $x' \prec x$\,. Suppose $x \in X_0$\,. But this contradicts the assumption that $x' \in N$. Hence, $x \in {\cal R}^n \setminus X_0$\,. \mybox

\vspace{0.2cm}
Lemma \ref{lemma7.2} shows how to select candidates for upper shells, which satisfy condition (\ref{eq7.3}) and condition (\ref{eq3.5}), or equivalently, condition (\ref{eq3.4}) and condition (\ref{eq3.5}).

The condition in Lemma \ref{lemma7.2} and Lemma \ref{lemma7.3} that all objective functions are strongly monotonically increasing cannot be relaxed, as illustrated in Figure 4. Dashed lines are contours of strongly monotonically increasing objective function $f_1$, dotted lines are contours of objective function $f_2$ which is not strongly monotonically increasing, the thick line shows set $N$. Function $f_1$ attains its maximum at $x^\star$ and function $f_2$ attains its maximum at $\bar{x}$, thus $x^\star \in N$ and $\bar{x} \in N$. Elements $x^1,\, x^2$ satisfy $\bar{x} \leq x$, however neither $x^1 \in {\cal R}^n \setminus X_0$ nor $\bar{x} \prec x^1$ holds, and $x^2 \in {\cal R}^n \setminus X_0$ holds but not $\bar{x} \prec x^2$. Figure 5 represents this situation in the space of objective function values. 
Analogous drawings can be made for any strongly monotonically increasing function and any set $X_0$.  

\begin{figure}[t]
\label{figure4}
\begin{center}
\rotatebox{0}{\includegraphics[height=2.0in]{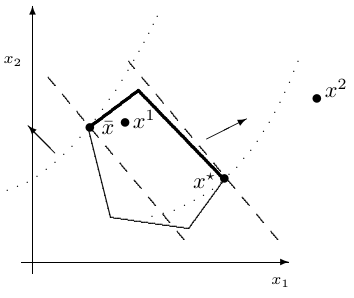}}
\end{center} 
\caption{An example illustrating why the assumption in Lemma \ref{lemma7.2} and in Lemma \ref{lemma7.3} that all objective functions are strongly monotonically increasing cannot be relaxed. }
\end{figure}

\begin{figure}[h]
\label{figure5}
\begin{center}
\rotatebox{0}{\includegraphics[height=2.0in]{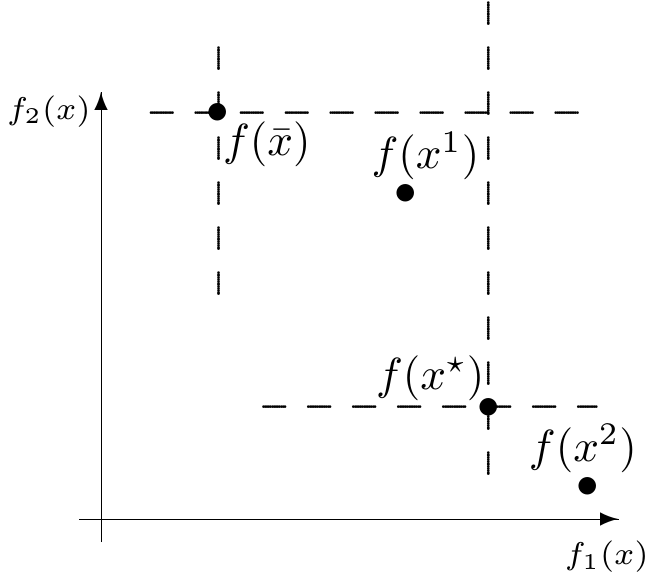}}
\end{center}
\caption{An example illustrating why the assumption in Lemma \ref{lemma7.2} and in Lemma \ref{lemma7.3} that all objective functions are strongly monotonically increasing cannot be relaxed -- the situation in the space of objective function values.}
\end{figure}

Lemma \ref{lemma7.3} is of the existential type since set $N$ is in general unknown. However, there is a class of MO problems, defined in Lemma \ref{theorem7.4}, in which $x' \in X_0, \ x \in {\cal R}^n \setminus X_0$ and $x' \leq x$ implies  $x' \prec x$. Thus, any subset of elements $x \in {\cal R}^n \setminus X_0$ which satisfy $g(x) > b$, $x' \leq x$ for some $x' \in X_0,$ and condition (\ref{eq3.3}), is a valid upper shell.
\begin{Theorem}
\label{theorem7.4}
Let all objective functions $f_l, \ l = 1,\dots,k$, be strongly monotonically increasing on ${\cal R}^n$\,.  Let one of the conditions defining $X_0$ be of the form
\[
g(x) \leq b \,,
\]
and let $g(x)$ be strongly monotonically increasing on ${\cal R}^n$. Then, $x' \in X_0$ and $x \in \{ x \ | \ g(x) > b \}$ and $x' \leq x$, implies $x' \prec x$.
\end{Theorem}

\noindent Proof. Any element $x$ such that $g(x) > b$ belongs to ${\cal R}^n \setminus X_0$. For any element $x'$ of $X_0$ and any element $x$ of ${\cal R}^n \setminus X_0$ such that $x' \leq x$, and for strongly monotonically increasing functions $f_l(x), \ l=1,\dots,k$, we have $f_l(x') < f_l(x)$. Elements $x$ such that $x \in {\cal R}^n \setminus X_0$ and $x' \leq x$ exist since $X_0$ is compact. Thus $x' \prec x$. \mybox

\vspace{0.2cm}
Theorem \ref{theorem7.4} is illustrated in Figure 6. Dashed and dotted lines are the contours of the strongly monotonically increasing functions $f_1$ and $f_2$, dashed thick line is the contour of the strongly monotonically increasing function $g(x)$.  For all infeasible elements $x \in \{x^1, x^2, x^3, x^4 \}$ relation $x' \prec x$ for some $x' \in X_0$ holds. This time no information on $N$ is assumed. Figure 7 represents this situation in the space of objective function values.

Theorem \ref{theorem7.4} is constructive. If the assumptions of the lemma hold, any subset of elements $x \in {\cal R}^n \setminus X_0$ which satisfy $g(x) > b$, $x' \leq x$ for some $x' \in X_0$, and condition (\ref{eq3.3}), is a valid upper shell.

\begin{figure}[t]
\label{figure6}
\begin{center}
\rotatebox{0}{\includegraphics[height=2.0in]{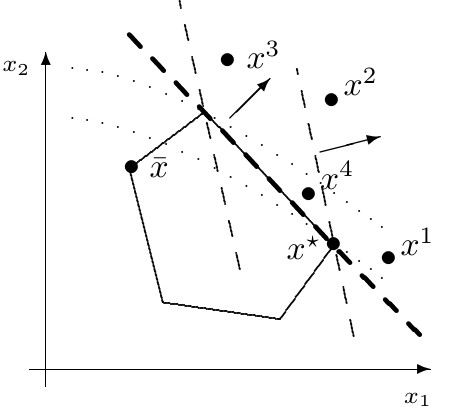}}
\end{center}
\caption{An illustration to Theorem \ref{theorem7.4}.}
\end{figure}
\begin{figure}[h]
\label{figure7}
\begin{center}
\rotatebox{0}{\includegraphics[height=2.0in]{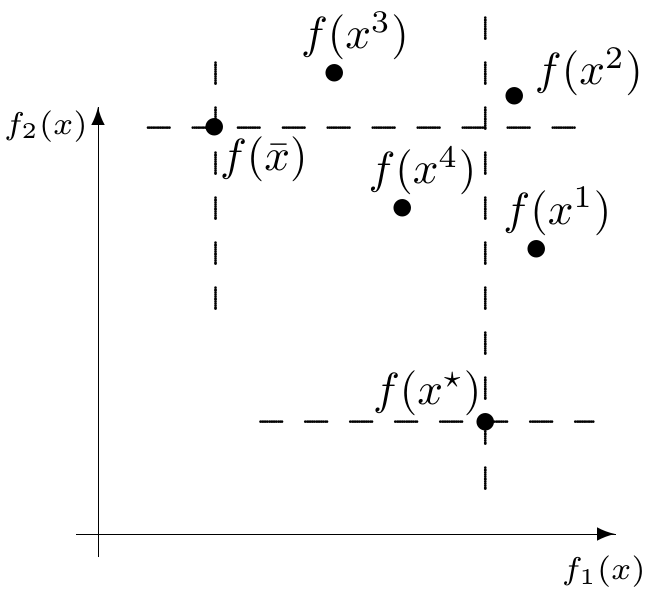}}
\end{center}
\caption{An illustration to Theorem \ref{theorem7.4} -- the situation in the space of objective function values.}
\end{figure}

Linear multiobjective and linear mixed-integer problems with positive coefficients in objective and constraint functions and $\leq$ type constraints are the simplest examples of problems which fall to this class. The property persists if linear functions are replaced by any strongly monotonically increasing function. As already mentioned above, this property has been exploited in the context of biobjective multidimensional knapsack problems \cite{Kaliszewski_2016} but the approach is directly extendable to any number of criteria.

\section{A Numerical Example}
\label{section8}

\vspace{0.2cm}
Consider the MO modeling problem -- a round beam with mass and deflection as objective functions (\cite{KKM_2016}).

\noindent
\textit{The MO problem.}

\vspace{0.2cm}









\[
\begin{array}{cc}
 &  ''max'' f(d,g) = \left( \begin{array}{l}
                  -f_1(d,g) = -\pi (d + g) g \rho l \\ \\
                  -f_2(d,g) = -\frac{4 F l^3}{3 E \pi ((d+2g)^4 - d^4)} 
                \end{array} \right) 
\end{array}
\]
\[
\begin{array}{cc}
 &      X_0 = \left\{ (d,g) \left| \begin{array}{c}
                             \frac{8 F l}{\pi} \frac{(d + 2g)}{(d + 2g)^4 - d^4} \leq k_g \\ \\
                             0 < d \leq 0.1  \\ \\ 0.001 \leq g \leq 0.1 
                            \end{array}  \right. \right.
\end{array}
\]
\noindent where

\begin{tabular}{ll}
$f_1(d,g)$ -- & (mass [kg]), \\

$f_2(d,g)$ -- & (deflection [m]), \\

$d$ -- &  (internal diameter [m]), \\

$g$ -- &  (wall thickness [m]), \\ 

$F = 10^4$ & (bending force [N]), \\ 

$l = 3$ & (beam length [m]), \\

$\rho = 7.86 \cdot 10^3$ & (material density [$\frac{\mbox{kg}}{\mbox{m}^3}$]), \\

$E = 2.1 \cdot 10^{11}$ & (Young modulus [Pa]), \\

$k_g = 150 \cdot 10^6$ & (maximal bending stress [Pa]).
\end{tabular}

\vspace{0.2cm}
To remain consistent with the problem formulation (\ref{eq3.1}) and the definition of the dominance relation, we maximize $- f_1(d,g)$ and $- f_2(d,g)$.

The second objective function is strongly monotonically increasing on ${\cal R}^2$ but the first is not. Thus, in this case an upper shell has been succesfully constructed by the combination of Lemma \ref{lemma7.1} and the relaxation approach (Lemma \ref{lemma5.1}).

It is also worth observing that by Lemma \ref{lemma6.1}, replacement of the second objective function by $ - f'_2(d,g) = -\frac{4 F l}{3 E \pi ((d+2g)^4 - d^4)}$ leaves the efficient set $N$ unchanged.

\begin{figure}[t]
\label{figure8}
\begin{center}
\rotatebox{0}{\includegraphics[height=2.5in]{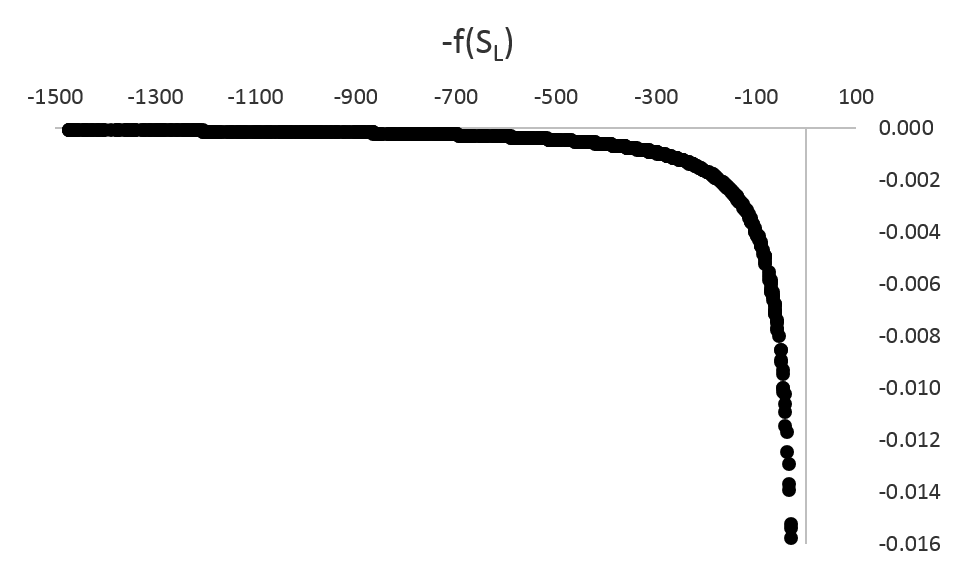}}
\end{center}
\caption{The image $-f(S_L)$ of a lower shell $S_L$ for the numerical example of Section \ref{section8} (mass -- horizontal axis).}
\end{figure}

\begin{figure}[h]
\label{figure9}
\begin{center}
\rotatebox{0}{\includegraphics[height=2.5in]{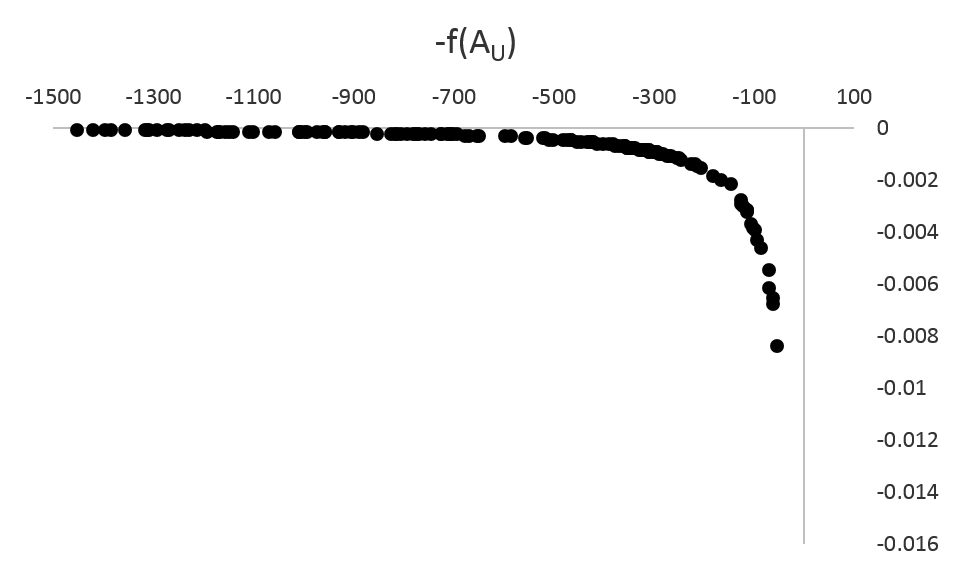}}
\end{center}
\caption{The image $-f(A_U)$ of an upper approximation $A_U$ for the numerical example of Section \ref{section8}  (mass -- horizontal axis). Here $A_U$ satisfies the definition of upper shell.}
\end{figure}

Figure 6 and Figure 7 present examples of a lower shell and an upper approximation, represented by the objective function mapping. They were derived by the algorithm described in \cite{KMP_2012}.

\section{Concluding Remarks}
\label{section9}

The results presented in the paper are inspired by attempts to provide tools for solving large and computationally expensive MO problems. As already said, with a pair of a lower shell $S_L$ and an upper shell $S_U$ it is possible to approximate selected efficient elements of $X_0$ with controllable accuracy (\cite{KM_2014,KMP_2012}).
In consequence, in the full analogy to singleobjective optimization, this enables stopping computations whenever satisfactory approximation accuracy is reached.
We have shown that with a rather mild conditions on problem (\ref{eq3.1}) there exist elements which form upper approximations to that problem.
We have also shown that some problem modifications, if admissible, guarantee that upper approximations have properties of upper shells. This adds to the fact that there are instances of problem (\ref{eq3.1}) where this is always the case. For example, multidimensional knapsack and multidimensional multiple choice knapsack problems, set covering (after a suitable transformation) and set packing problems have this feature if infeasible $x$ are confined to $\{ 0, 1 \}^n$.  Identification of other classes of problems for which upper approximations have properties of upper shells will be the subject of our further research.

Lemma \ref{lemma6.1} can be particularly useful in large-scale computations, where the cost of computing $f(x)$ becomes a limiting factor. Lower shells and upper approximations can be derived with functions $f'_l(x)$ which generate the same linear order as functions $f_l(x)$, but of lower computing cost. Linear functions and polynomial functions when defined on appropriate domains can serve here as the simplest example.

One might rightly argue that a natural vehicle to implement the concept of the two-sided Pareto front approximations is evolutionary multiobjective optimization, as presented e.g. in the monographs \cite{Deb_2001,Coello_2002},
and numerous papers published on the subject. Moreover, in our earlier papers (\cite{KM_2014,KMP_2012,KM_2010,KMP_2011,KKM_2016})
we have made use of this vehicle. However, seeing EMO as a natural but not necessarily the only mechanism to populate PF approximations, in this work we purposely have not related directly our results to that specific kind of heuristics, because they are applicable to heuristics (to derive lower and upper shells) of any sort.

\end{document}